\documentclass[12pt]{article}
\usepackage{amsmath}
\usepackage[dvipsnames,usenames]{color}
\usepackage{amsfonts}
\usepackage{mathrsfs}
\usepackage{amssymb}
\usepackage{amsthm}
\def\endpf{\hbox{\vrule height1.5ex width.5em}}

\def\<{\langle}
\def\>{\rangle}
\numberwithin{equation}{section}
\def\<{\langle}
\def\>{\rangle}

\def\p{\partial}
\def\a{\alpha}

\def\O{\Omega}

\def \sm{\setminus}
\def\-{\overline}

\def\e{\epsilon}

\def\endpf{\hbox{\vrule height1.5ex width.5em}}

\def\a{\alpha}
\def\endpf{\hbox{\vrule height1.5ex width.5em}}

\def\-{\overline}

\def\O{\Omega}

\def\sm{\setminus}

\def\wt{\widetilde}

\def\endpf{\hbox{\vrule height1.5ex width.5em}}

\def\a{\alpha}

\def\a{\alpha}

\def\endpf{\hbox{\vrule height1.5ex width.5em}}
\def\a{\alpha}

\newtheorem{theorem}{Theorem}[section]
\newtheorem{lemma}[theorem]{Lemma}
\newtheorem{Corollary}[theorem]{Corollary}
\newtheorem{proposition}[theorem]{Proposition}

\newtheorem{Definition}[theorem]{Definition}

\newtheorem{remark}[theorem]{Remark}

\date{\ }

\begin{document}
\title{\bf $\overline\partial$-equation on a lunar domain  with mixed boundary conditions}

\author{{Xiaojun  Huang}\footnote{
Supported in part by NSF-1101481}\qquad  Xiaoshan
Li\footnote{Supported by the China Scholarship Council and the
Fundamental Research Fund for the Central Universities}}

\maketitle

\begin{abstract}\vskip 3mm\footnotesize
\noindent In this paper,  making use of the method developed by
Catlin,
 we study the $L^2$-estimate for the $\bar\partial$-equation
 on a lunar manifold with the mixed boundary conditions.

\vskip 4.5mm



\end{abstract}

\vskip 12mm

\bigskip
\tableofcontents
\section{Introduction}

In this paper, we study the $\-{\p}$-equation for $(0,q)$-forms  on
a special type of non-smooth domain $S^{\e_0}_\varphi$, called a lunar domain,
in a complex manifold with mixed boundary conditions. The domain we
are considering here has two pieces  of the boundaries $M_0$ and $M_1$
intersecting highly tangentially along a smooth real-submanifold
$E$. We assume that $M_0$ has at least $(q+1)$-positive Levi
eigenvalues or $(n-q+1)$-negative Levi eigenvalues. Assume $M_1$ has
the opposite property for the Levi eigenvalues as that for $M_0$.

We impose the $\overline\partial$-Dirichlet boundary condition on $M_0$  and the
$\-{\p}$-Neumann boundary condition on $M_1$. We introduce a
Hermitian metric over $S^{\e_0}_\varphi$ such that $E$ can be treated  as the
infinity of $S^{\e_0}_\varphi$. We will establish an $L^2$-estimate and derive a
Hodge-type decomposition theorem in this setting.

$\-\p$-equations over such a special type of non-smooth domains,
with mixed boundary conditions, are of fundamental importance in
understanding many geometric problems. In the deep papers of Catlin
[Cat], Cho [Cho] and Catlin-Cho [CC], such equations played a
crucial role for studying various extension problems for CR
structures, which are directly linked to the local embedding problem
of abstract CR manifolds with certain signature conditions. In a
paper of Huang-Luk-Yau [HLY], solving such a $\-\p$-equation for
$(0,2)$-forms  also played an important role for the study of
various deformation problems for compact strongly pseudoconvex CR
manifolds of at least five dimension. In the work  of Catlin [Cat],
Catlin-Cho [CC] and Cho [Cho], the domain encountered  is only
assumed to be sitting in an almost complex manifold. However,  the
domain is  uniformly scaled such that it is  sufficiently close to
$M_0$. In this specific setting, Catlin proved that there is no
cohomology obstruction for solving the $\-\p$-equations. Other
related studies for the $\bar{\p}$-Dirichlet problem can be found in
the work of Chakrabarti-Shaw [CSh].

In this paper, we will study the above mentioned $\-\p$-equation,
with the mixed boundary conditions, without any scaling of the
domain. Then one does not expect the $\-\p$-equation is always
solvable. However, we will show that the obstruction is of finite
dimension. Though we basically follow the approach of Catlin [Cat],
one key point in our paper is that we use the property close to the
non-smooth corner near $E$ for a different weighted metric to avoid
the difficulty which was circumvented  in [Cat], [Cho] only by
uniformly shrinking the lunar domain $S^{\e_0}_\varphi$ toward
$M_0$.

$\-{\p}$-equations with various boundary conditions are the basic
tools to work on many geometric or analytic problems in Several
Complex Variables and Complex Geometry. There is a vast amount of
work done in the literature. Here, we only refer the reader to  the
books by Folland-Kohn [FK], H\"ormander [Ho2], Demailly [DE] and Chen-Shaw [CS],
as well as, many references therein.

\section{Basic set-up and statement of the main theorem}
Let $M$ be a smooth hypersurface of real dimension
$2n-1(n\geq3)$ in a  complex manifold $X$ of real dimension $2n$. Let
$\varphi\in C^\infty(M)$ be a function such that $d\varphi(x)\neq0$
when $\varphi(x)=0$.
 Write $K=\overline{\{x\in M: \varphi(x)> 0\}}$. Assume $K\subset\subset M$ is bounded domain in $M$
 with smooth boundary  $E=\{x\in M: \varphi(x)=0\}$.

For a sufficiently small
 $\e_0>0$, define $M^{\e_0}=\{x\in M: |\varphi(x)|<\e_0\}.$ Suppose
that there exists a tubular neighborhood $\mathcal{N}_{\e_0}$ of
$M^{\e_0}$ in $X$ and a $C^\infty$ map $\Phi$ such that
$\Phi:\mathcal{N}_{\e_0}\rightarrow M^{\e_0}\times(-2,2)$ is a
diffeomorphism. Write $\Omega_{\e_0}=M^{\e_0}\times (-2,2)$, $
\mathcal{L}=\Phi_\ast(T^{1,0}\mathcal{N}_{\e_0})$, where
$T^{1,0}\mathcal{N}_{\e_0}$ is the holomorphic tangent bundle of
$\mathcal{N}_{\e_0}$. Then $(\Omega_{\e_0},\mathcal{L})$ is a
complex manifold biholomorphic to
$(\mathcal{N}_{\e_0},T^{1,0}\mathcal{N}_{\e_0})$. Also
$(M^{\e_0}\times \{0\},\mathcal{L}|_{(M^{\e_0}}\times\{0\})\cap
\mathbb{C}T(M_0\times \{0\})$ is a CR hypersurface  in
$(\Omega^{\e_0}, \mathcal{L})$. Identify $M^{\e_0}\times\{0\}$ with
$M^{\e_0}$.  Write ${
S}=\mathcal{L}|_{M^{\e_0}}\cap\mathbb{C}TM^{\e_0}$ which is the $CR$
bundle of $M^{\e_0}$. In what follows, when there is no risk of
causing confusion, we identify ${\cal N}_{\e_0}$ with $\O_{\e_0}$
and objects defined over ${\cal N}_{\e_0}$ with those corresponding
ones over $\O_{\e_0}$.

Define $r(x,t)=t\varphi^{-4}(x)$. Assume $S^{\epsilon_0}_\varphi$ is
a bounded domain in $X$ with two pieces of connected boundaries
$M_0:=M\cap \{\varphi> 0\}$ and $M_1$, whose closures intersect $M$
tangentially  along $E$. Moreover,
 $S_\varphi^{\epsilon_0}\cap {\mathcal N}_{\epsilon_0}: =\{(x,t)\in
\Omega_{\epsilon_0}|\varphi(x)>0, -1< r(x,t)< 0\}$ and
$M_1\cap {\mathcal N}_{\e_0} =\{(x,t)|\varphi(x)>0,r(x,t)=-1\}$.

Equip $X$ with a Hermitian metric. For any $x_0\in M_0$ or $x_0\in
M_1$, let $\{L_j\}_{j=1}^n$ be a smooth orthornormal basis of the
cross sections of $T^{(1,0)}(W(x_0))$, where $W(x_0)$ is
sufficiently small neighborhood of $x_0$ in the ambient space. Let
$\{\omega_j\}_{j=1}^n$ be its dual frame. Assume that $L_j$ are
tangent to $M$ or $M_1$, when restricted to $M_0$ or $M_1$, for
$j\neq n$, respectively. For a $(0,q)$-form with $0<q\le n$
$$U=\sum_{j_1<j_2<\cdots<j_q}U^{j_1j_2\cdots
j_q}\overline{\omega_{j_1}}\wedge\cdots\wedge\overline{\omega_{j_q}}$$
defined in the side of $W(x_0)\cap M_0$ or of $W(x_0)\cap M_1$,
which is  inside $S_\varphi^{\e_0}$, that is smooth up to $M_0$ or
$M_1$. We say $U$ satisfies the $\overline\partial$-Dirichlet
condition along $M_0$ if $U_J|_{M_0}\equiv0$ whenever
$J=(j_1,\cdots,j_q)$ with $j_q\neq n$. We say $U$ satisfies the
$\overline\partial$-Neumann condition along $M_1$ if
$U_J|_{M_1}\equiv0$ when $j_q=n$. Apparently, this definition is
independent of the choice of the Hermitian metric over $X$. Indeed,
one only needs a smooth Hermitian metric over
$\-{S_\varphi^{\e_0}}\sm E$ to define the Dirichlet or Neumann
boundary conditions along $M_0$ and $M_1$.

 Following Catlin in [Cat], we write $\mathcal{E}_c^{(0,q)}$ for the collection of smooth $(0,q)$-forms with compact support in
  $\-{S^{\e_0}_\varphi}\sm E$.
 Write $\mathcal{B}_+^{q}(S^{\e_0}_\varphi)$ for the subset of $\mathcal{E}_c^{(0,q)},$ whose elements satisfy
 $\overline\partial$-Dirichlet boundary condition along $M_0$. Write $\mathcal{B}_-^q(S^{\e_0}_\varphi)$
for the subset of  $\mathcal{E}^{(0,q)}$ whose elements satisfy the $\overline\partial$-Neumann boundary condition along $M_1$.
Now, we will use the specific Hermitian metric over ${S}^{\e_0}_\varphi$ to be defined in (\ref{eqn6}) of Section
3, which is smooth up to the boundary $M_0\cup M_1\sm E$  and blows
up at a suitable rate
when approaching their intersection $E$.
We define $L_{(0,q)}^2(S_\varphi^{\e_0})$ to be the space of $(0,q)$-forms with coefficients being $L^2$-integrable with respect to this metric. We extend the $\overline\partial$-operator to the $L^2$-space in
the following way:

We say that $U\in L_{(0,q)}^2(S_\varphi^{\e_0})$ is in the domain of
the operator $T$ with $TU=F$ if for any
$V\in\mathcal{B}_{-}^{q+1}(S_\varphi^{\e_0})$, we have
$(U,\overline\partial^\prime V)=(F,V)$, where
$\overline\partial^\prime$ is the standard formal adjoint operator
of $\overline\partial$ with respect to this specific Hermitian
metric. Similarly, we define $S:\
L_{(0,q-1)}^2(S_\varphi^{\e_0})\rightarrow
L_{(0,q)}^2(S_\varphi^{\e_0})$ and let $T^\ast$ and $S^\ast$ be
their Hilbert adjoints. We define $Q(U,U)=\|TU\|^2+\|S^{\ast}U\|^2$
to  be the $Q$-norm associated with the operators $T$ and $S^\ast$.

In [Cat] and [Cho], to study  the extension of CR structure of $M$,
the authors obtained a standard
$L^2$-estimate with respect to the $\overline\partial$-operator with
mixed boundary condition when the thickness of $S_{\varphi}^{\e_0}$   is
sufficiently small. (See [Corollary 7.10, Cat]).

In this paper, we consider the $L^2$-estimate with respect to a
$\overline\partial$-operator with mixed boundary conditions.
However, the thickness of $S^{\e_0}_\varphi$ can be arbitrary.
Define $N^+(K)$ (respectively,  $N^-(K)$) to be the largest $m\geq
0$ such that the Levi form has at least $m$ positive (respectively,
negative) eigenvalues at each $x\in K$ with respect to the domain
${S_\varphi^{\e_0}}$. Define $N^+(M_1)$ (respectively, $N^{-}(M_1)$)
has the same meaning as for $N^{+}(K)$ (respectively, for
$N^{-}(K)$) with respect to ${S_\varphi^{\e_0}}$, too. Then our
main theorem is the following:

\begin{theorem}\label{101}
Assume the above notations and definitions. If either we have
$N^+(K)\geq q+1$  and $N^-(M_1)\geq q+1$ or we have $N^-(K)\geq
n-q+1$ and $N^+(M_1)\geq n-q+1$. Then there exists a neighborhood
$V_{c,0,1}$ of the boundary of $S_\varphi^{\e_0}$ in
$\-{S_\varphi^{\e_0}}$ and a constant $C>0$ such that for any $U\in
L^2_{(0,q)}(S^{\e_0}_\varphi)$ with $U\in Dom(T)\cap Dom(S^{\ast})$,
it holds that
\begin{equation}\label{eqnc}
\int_{V_{c,0,1}}|U|^2dV\leq C\left (
Q(U,U)+\int_F|U|^2dV\right ),
\end{equation}
where F is a certain  compact subset of $S^{\e_0}_\varphi$
independent of $U$.
\end{theorem}

\begin{Corollary}\label {202}
Write $H^{(0,q)}(S_\varphi^{\e_0})$ for  the quotient space
$N_T/R_S$ with $ N_T=\{U:U\in L_{(0, q)}^2(S_\varphi^{\e_0}),\
TU=0\}$ and $R_S$ the image of the operator $S$. Then
$H^{(0,q)}(S_\varphi^{\e_0})=N_T/R_S$ is of finite dimension.
\end{Corollary}

\section{Existence of the special frames on $S_\varphi^{\e_0}$ near $E$}
For the proof of Theorem \ref{101}, we follow the approach in Catlin
[Cat] and Catlin-Cho [CC]. However, we need to choose a different
weight of blowing up for the metric near the singular set $E$ of the
boundary to deal with the difficulty caused by not shrinking the
thickness of the lunar domain. This also requires the modification
for the choice of  the special frame to study the $L^2$-estimates
later. For the convenience of the reader, we give a detailed exposition
on the choice of the frame in this section.

Let $M_{t_0}$ near $E$  be defined by the defining equation $t=t_0$.
Define $\eta=-\frac{1}{2}(i \p t-i\-{\p}t)$ near $E$. Then $\eta$ is
a real-valued 1-form and is a contact form along each $M_t$ near
$E$. Let $X_0$ be a real-valued smooth vector field tangent to $M$
near $E$ such that $(\eta,X_0)=1$ over $M$ near $E$. Extend $X_0$ to
a neighborhood of $E$ in $\O_{\e_0}$, independent of $t$, and scale $X_0$
if needed. Then we can get a real-valued smooth vector field $X_0$
in a neighborhood of  $E$ in $\O$ such that  near $E$ $(\eta,X_0)=1$
and $X_0(t)\equiv 0$.

We assume, without loss of generality, that the Levi form of $M_0$
is defined by $\sqrt{-1}\eta([X_1,\overline{X_2}]),$ $ X_1, X_2\in
S.$
Write $S_{(x,t)}$ for the subspace of $\mathcal{L}_{(x,t)}$ that are
tangent to $M_t$ near $E$.
Set $Y_0=-J_{\mathcal{L}}(X_0)$, so that $X_0+\sqrt{-1}Y_0$ is a
section of $\mathcal{L}$ that is transversal to the level set $t$.
Let $G:\Omega_{\e_0}\cap O(E)\rightarrow\Omega_{\e_0}\cap O(E)$ be a
diffeomorphism such that $G$ fixes $M\cap O(E)$ and
$$G_{\ast}Y_0|_{(x,0)}=\frac{\partial}{\partial t}\Big|_{(x,0)},\ \ x\in
M\cap O(E).$$ Here we write $O(E)$ for a small neighborhood of $E$
in $\O_{\e_0}$. Since $dt(J_{\mathcal{L}}(X_0))$ always has the same
sign (If not, $X_0+\sqrt{-1}Y_0$ is a section that is tangent to the
level set), we may assume that $dt(J_{\mathcal L}(X_0))<0$, thus
$dt(Y_0)>0$ along $M_0$. Hence  $G$ preserves the sides of $M_0$.
Then $\widetilde{Z}=-\sqrt{-1}G_{\ast}(X_0+\sqrt{-1}Y_0)$ is a
global section of $S_\varphi^{\e_0}$ near $E$ such that along $M_0$
\begin{equation}\label{eqn1}
\widetilde{Z}=-\sqrt{-1}X_0+\frac{\partial}{\partial t}.
\end{equation}
We write  $\widetilde Z=\wt{X}+g(x,t)\frac{\partial}{\partial t},$ where
$\wt{X}t\equiv0.$ Then we set $ Z_n=X+\frac{\partial}{\partial t}$
near $E$ with $X=g^{-1}(x,t)\wt{X}$.

We define another subbundle of
$\mathcal{L}$ on $S_{\varphi}^{\e_0}$ by setting
\begin{equation}\label{eqn4}
\mathcal{R}_{(x,t)}=\{L\in\mathcal{L}_{(x,t)}: Lr=0,\
r=t\varphi^{-4}(x)\}.\
\end{equation}
Clearly, the map defined by
\begin{equation}\label{eqn5}
H(L)=L-L(r)(Z_nr)^{-1}Z_n,\  L\in S_{(x,t)},
\end{equation}
defines an isomorphism of $\mathcal{S}:=\cup_{(x,t)\approx E}
S_{(x,t)}$ onto $\mathcal{R}:=\cup_{(x,t)\approx E} R_{(x,t)}$, where
$$Z_n(r)=\varphi(x)^{-4}\left(1+(-4t)\varphi(x)^{-1}X\varphi(x)\right).$$

We fix a smooth Hermitian metric $\langle,\rangle_0$ on
$\overline{S_\varphi^{\epsilon_0}}$ that is induced from the
Hermitian metric on $X$ such that on $\Omega_{\e_0}$ we have
$\langle Z_n,Z_n\rangle_0=1$ near $E$. We define a new Hermitian metric $\langle,\rangle$ on
$\-{{ S}_\varphi^{\e_0}}\sm E$ such that near $E$ we have the
following relations :

\begin{equation}\label{eqn6}
\begin{split}
\langle H(L_1),H(L_2)\rangle&=\varphi^{-4+\lambda}(x)\langle L_1,L_2\rangle_0, L_1,L_2\in \mathcal{S}\\
\langle H(L_1), Z_n\rangle&=0, L_1\in \mathcal{S}\\
\langle Z_n,Z_n\rangle&=\varphi^{-8+2\lambda}(x)
\end{split}
\end{equation}
where $\lambda$ is a constant with $0<\lambda<\frac12.$ We now show
that  $\overline {S_\varphi^{\epsilon_0}}\setminus E$ near $E$ can be covered by special
coordinate systems such that on each chart there is an orthonormal
frame of $\mathcal{L}$ that satisfies good estimates. This is
fundamentally important for it then helps to reduce the non-compact
situation to more or less the compact situation. Comparing with the
weight in [Cat], we add $\varphi^\lambda(x)$ to take care of the
trouble created from the corner near $E$.

\begin{proposition}\label{thm1}
For any $x_0\in M_0$ with $0<\varphi(x_0)\ll 1$, there exists a
neighborhood $W(x_0)\subset \overline {S_\varphi^{\e_0}}\setminus E$
with the following properties:

(i) On $W(x_0)$, there are smooth coordinates $y_1,\ldots,y_{2n}$
so that
\begin{equation}\label{eqn7}
W(x_0)=\left\{y:|y^\prime|<\sigma_0, -\varphi^\lambda(x_0)\leq
y_{2n}\leq 0\right\},
\end{equation}
where $\sigma_0$ is a constant independent of $x_0$ to be determined
later. Also  $y^\prime=(y_1,\ldots,y_{2n-1})$ is independent of $t$
and $y_{2n}=t\varphi^{-4}(x)\varphi^{\lambda}(x_0).$ $M_0\cap
W(x_0)$ and $M_1\cap W(x_0)$ correspond to  points in $W(x_0)$ with
$y_{2n}=0$ and $y_{2n}=-\varphi^\lambda(x_0),$ respectively.
Moreover, the point $x_0$ corresponds to the origin.

(ii) On $W(x_0)$, there exists a smooth orthonormal  frame
$L_1,\cdots,L_n$ for $\mathcal{L}$ such that if $\omega^1,\cdots,
\omega^n$ are the dual frame, and if $L_k$ and $\omega^k$ are written
as $\sum\limits_{j=1}^{2n}b_{kj}\frac{\partial}{\partial y_j}$ and
$\sum\limits_{j=1}^{2n}d_{kj}dy_j,$ respectively, then

\begin{equation}\label{eqn9}
\sup\limits_{y\in W(x_0)}\left\{|D_y^\alpha b_{kj}(y)|+|D_y^\alpha
d_{kj}(y)|\right\}\leq C_{|\alpha|},\ \  b_{k 2n}=0\ \hbox{for}\
k\not = n.
\end{equation}
where $C_{|\alpha|} $ is independent of $x_0,j,k$.

(iii) If $d_1\leq d_2\leq\cdots\leq d_{n-1}$ are the eigenvalues
of the Levi form $\sqrt{-1}\eta([L_1,\overline{L_2}])$ at $x_0$,
then at every point $x\in W(x_0)$, we have the following estimates
\begin{equation}\label{eqn10}
\begin{split}
\omega^n([L_i,\overline{L_j}])(y)&=O(\sigma_0),
i\neq j, i,j<n,\\
\omega^n([L_j,\overline{L_j}])(y)&=d_j+O(\sigma_0), j<n.
\end{split}
\end{equation}
Here we use $O(\sigma_0)$ to denote the terms which are bounded by $C\sigma_0$ where $C$ is a constant independent of $x_0$ and $\sigma_0$.

(iv) For each sufficiently small $\sigma_0$, there is a countable
family $\{W(x_\a)\}$ such that it covers a fixed neighborhood of
$K\sm E$ in $K$ near $E$ and for any point $p\in K\sm E$ in this
neighborhood, there are at most $N_0$ elements from this family that
contain $p$. Here $N_0$ is independent of the choice of  $p$ and
$\sigma_0$. Moreover, for each $\alpha$, there is a function
$\xi_\alpha\in C^{\infty}_0(W(x_\alpha)\cap M)$ such that
$\sum_{\a}\xi_\a^2\equiv 1$ and the differentiation of $\xi_\a$ with
respect to the $y'$-coordinates is bounded by $C/\sigma_0$ with $C$
a fixed constant independent of $x_\a$ and $\sigma_0$.
\end{proposition}
\begin{proof}
The proof of this proposition is similar to that in [Cat], though
adding a new scale $\varphi(x_0)^\lambda$ requires  modifications.
For convenience of the reader, we include  all the  details. First,
there exists a finite number of coordinate charts $V'_v,
v=1,\cdots,N$ in $M$ that cover K near $E$ in $M$ such that on each
$V'_v$, there exists coordinates $(x_1,\ldots,x_{2n-1})$ with
$\frac{\partial}{\partial x_{2n-1}}=-X_0$ at all points in $V'_v$.
Also, $V'_v$ is defined by $|x'|<\e_1$ for a certain fixed small
$\e_1>0$. Define $V_v:=V'_v\times [0,-1]$ and set on $V_v$,
$x_{2n}=t$, $x_k(x^\prime,t)=x_k(x^\prime),$  $k<2n$ for
$x^\prime\in V'_v$. We can assume  that there exists an orthonormal
frame $\{L_i^v\}_{i=1}^{n-1}$ of $\mathcal{S}$ with respect to the
former fixed Hermitian metric $<,>_0$ of $\mathcal{L}$ in $V_v$. Let
$L_n^v=Z_n$. For any point $x_0\in M$ with $0<\varphi(x_0)<<1$, by
the Lesbeque covering lemma, we can assume that $x_0\in V_v$ for  a
certain $v$ with $|x'(x_0)|<\e_2$, where $0<\e_2<\e_1$ is
independent of $x_0$.
 We
can define an affine transformation
$C_{x_0}^v:\mathbb{R}^{2n}\rightarrow\mathbb{R}^{2n}$ so that if
$(x_0^\prime,0)\in \mathbb{R}^{2n}$ is the coordinates of $x_0$,
then
\begin{equation}\label{eqn11}
C_{x_0}^v(x^\prime, x_{2n})=(P_{x_0}(x^\prime-x_0^\prime),x_{2n}),
\end{equation}
where $P_{x_0}$ is a $(2n-1)\times(2n-1)$ constant matrix such that
in the new coordinates $\tilde x=(\tilde x_1,\cdots,\tilde x_{2n})$,
we have
\begin{equation}\label{eqn12}
\begin{split}
L_k^v|_{x_0}&=\frac{\partial}{\partial\tilde
x_{2k-1}}\Big|{x_0}-\sqrt{-1}
\frac{\partial}{\partial\tilde x_{2k}}\Big|_{x_0}, (1\leq k\leq n-1),\\
X_0\Big|_{x_0}&=-\frac{\partial}{\partial\tilde x_{2n-1}}\Big|_{x_0}.
\end{split}
\end{equation}

\noindent Also the domain where $\wt{x'}$ is defined contains a
fixed ball centered at the origin for any choice of $x_0$. Notice
that the second equality in (\ref{eqn12}) implies that
$X_0|_{(x^\prime,0)}=-\frac{\partial}{\partial \tilde
x_{2n-1}}|_{(x^\prime,0)}$ at all points of $M\cap V_v$. Hence,
along $M\cap V_v$,
\begin{equation}\label{eqn13}
L_n^v\mid_{(x^\prime,0)}=\sqrt{-1}\frac{\partial}{\partial\tilde
x_{2n-1}} \Big|_{(x^\prime,0)}+\frac{\partial}{\partial\tilde
x_{2n}}\Big|_{(x^\prime,0)}.
\end{equation}

\noindent We now define a new coordinates $y=(y_1,\cdots, y_{2n})$
by means of a dilation map $D_{x_0}:\mathbb{R}^{2n}\rightarrow
\mathbb{R}^{2n}.$ Set
\begin{equation}\label{eqn14}
\begin{split}
y&=D_{x_0}(\tilde{x})\\
&=\left(\varphi^{-2+\frac{\lambda}{2}}(x_0)\tilde x_1,\cdots,
\varphi^{-2+\frac{\lambda}{2}}(x_0)\tilde
x_{2n-2},\varphi^{-4+\lambda}(x_0) \tilde
x_{2n-1},\varphi^\lambda(x_0)\varphi^{-4}(x)\tilde x_{2n}\right).
\end{split}
\end{equation}

\noindent In terms of the $y$-coordinates, we define an open set
$W(x_0)$ by
\begin{equation}\label{15}
W(x_0)=\left\{x\in V_v\cap S_\varphi^{\e_0}:|y_k(x)| <\sigma_0,
k=1,\dots,2n-1, -\varphi^\lambda(x_0) \leq y_{2n}\leq 0\right\}.
\end{equation}
When $0<\varphi(x_0)<<1$, one can  apparently  find a fixed small
number $\sigma'_0>0$ such that every $W(x_0)$ is contained in some $V_v$
whenever $\sigma_0<\sigma_0'$. Notice that in $W(x_0)$, the set
where $y_{2n}=0$ and $y_{2n}=-\varphi^\lambda(x_0)$ coincides with
the set where $r(x,t)=0$ and $r(x,t)=-1$, respectively, which
represents the two boundaries of $S_\varphi^{\e_0}$.

Define a frame $L_1,\cdots, L_n$ on $W(x_0)$ by setting
\begin{equation}\label{eqn16}
\begin{split}
L_k&=\varphi^{2-\frac{\lambda}{2}}(x)\left(L_k^v-r_k(x)L_n^v\right)
=\varphi^{2-\frac{\lambda}{2}}(x)
H(L_k^v), k<n,\\
L_n&=\varphi^{4-\lambda}(x)L_n^v,
\end{split}
\end{equation}
where
\begin{equation}
r_k=(L_k^vr)(L_n^vr)^{-1}.
\end{equation}
Then $\{L_k\}_{k=1}^n$ forms  an orthonormal frame on $W(x_0)$ with
respect to the scaled Hermitian metric, and $\{L_k\}_{k=1}^{n-1}$
forms an orthonormal basis for $\mathcal{R}$.

If we write $L_k^v$ in terms of the $\tilde x$-coordinates
corresponding to $x_0$ as
\begin{equation}\label{eqn17}
\begin{split}
L_k^v&=\sum_{l=1}^{2n-1}e_{kl}(\tilde x)
\frac{\partial}{\partial \tilde x_l}, k<n,\\
L_n^v&=\frac{\partial}{\partial\tilde x_{2n}}
+\sum_{l=1}^{2n-1}e_{nl}(\tilde x)\frac{\partial}{\partial\tilde
x_l},
\end{split}
\end{equation}
and if we set
\begin{equation}\label{eqn18}
\begin{split}
E_{kl}&=e_{kl}\circ D_{x_0}^{-1}(y),R_k=r_k\circ D_{x_0}^{-1}(y), \varphi_l=\frac{\partial \varphi}{\partial\tilde x_l}\\
\Phi&=\varphi\circ D_{x_0}^{-1}(y), \Phi_l=\varphi_l\circ D_{x_0}^{-1}(y),
\end{split}
\end{equation}
then a calculation shows that
\begin{equation}\label{eqn19}
L_k^v-r_kL_n^v=\sum_{l=1}^{2n-1}\left(e_{kl}(\tilde{x})-r_k(\tilde{x})
e_{nl}(\tilde{x})\right)\frac{\partial}{\partial\tilde
x_l}-r_k(\tilde{x}) \frac{\partial}{\partial\tilde x_{2n}}
\end{equation}
and that the Jacobian matrix $\text Jac(D_{x_0})$ of $D_{x_0}$ is
\begin{equation}\label{20}
\begin{split}
\left[
\begin{array}{ccccc}
 \varphi^{-2+\frac{\lambda}{2}}(x_0) & 0& \cdots & 0 & 0 \\
 \vdots & \ddots& \vdots & \vdots & \vdots \\
 \vdots& \vdots & \varphi^{-2+\frac{\lambda}{2}}(x_0) & 0 & 0 \\
 0 & \cdots & 0 & \varphi^{-4+\lambda}(x_0) & 0\\
 -4\varphi_1\varphi^{-5}\varphi^\lambda(x_0)\tilde x_{2n} & \cdots& -4\varphi_{2n-2}\varphi^{-5}\varphi^\lambda(x_0)\tilde x_{2n} & -4\varphi_{2n-1}\varphi^{-5}\varphi^\lambda(x_0)\tilde x_{2n}& \varphi^\lambda(x_0)\varphi^{-4}
 \end{array}
 \right]
 \end{split}
\end{equation}
\noindent We conclude that in the $y$-coordinates of $W(x_0)$ when
$1\leq k\leq n-1,$
\begin{equation}\label{eqn21}
\begin{split}
&L_k=\sum_{l=1}^{2n-2}\frac{\Phi^{2-\frac{\lambda}{2}}(y)}
{\varphi^{2-\frac{\lambda}{2}}(x_0)}\left(E_{kl}(y)-R_kE_{nl}(y)\right)
\frac{\partial}{\partial y_l}\\
&+\frac{\Phi^{2-\frac{\lambda}{2}}(y)}
{\varphi^{4-\lambda}(x_0)}(E_{k,2n-1}(y)-R_kE_{n,2n-1}(y))\frac{\partial}
{\partial y_{2n-1}},\\
&L_n=\sum_{l=1}^{2n-2}\frac{\Phi^{4-\lambda}(y)}
{\varphi^{2-\frac{\lambda}{2}}(x_0)}E_{nl}(y)\frac{\partial}
{\partial y_l}+\frac{\Phi^{4-\lambda}(y)}
{\varphi^{4-\lambda}(x_0)}E_{n,2n-1}(y)\frac{\partial}
{\partial y_{2n-1}}\\
&+\left(\frac{\varphi^\lambda(x_0)}{\Phi^\lambda(y)}-4\sum_{l=1}^{2n-1}
\Phi^{3-\lambda}
\Phi_ly_{2n}E_{nl}(y)\right)\frac{\partial}{\partial y_{2n}}.
\end{split}
\end{equation}
Here
\begin{equation}\label{eqn22}
R_k=\frac{\sum\limits_{l=1}^{2n-1}(-4)E_{kl}(y)\frac{\Phi^3}
{\varphi^\lambda(x_0)}\Phi_l(y)y_{2n}}
{1+\sum\limits_{l=1}^{2n-1}(-4)E_{nl}(y)\frac{\Phi^3}
{\varphi^\lambda(x_0)}\Phi_l(y)y_{2n}}.
\end{equation}
Observe that  the diameter in the $\tilde{x}$-coordinates of
$W(x_0)$ is of the
quantity: $O(\varphi^{2-\frac{\lambda}{2}}(x_0))\ll\varphi(x_0)$ when
$0<\varphi(x_0)\ll1$.

Let $f\in C^{\infty}(W(x_0))$, we define
\begin{equation}
|f|_{m,W(x_0)}=\sup\limits_{y\in W(x_0)}\{|D_y^\alpha
f(y)|:|\alpha|\leq m\}
\end{equation}
and we can extend this norm to vector fields and $1$-forms by using
coefficients of $\frac{\partial}{\partial y_j}$ or $dy_j$. It can be
easily verified that
\begin{equation}
\lim\limits_{x_0\rightarrow E}|E_{kl}-b_{kl}|_{m,W(x_0)}=0,
\end{equation}
where $(b_{kl})_{n\times 2n}$ is a constant matrix given by
\begin{equation}\label{eqn23}
\begin{split}
&b_{k,2k-1}=1,b_{k,2k}=-\sqrt{-1}, k=1,\ldots,n-1,\\
&b_{n,2n-1}=\sqrt{-1},b_{n,2n}=1,
\end{split}
\end{equation}
and $b_{kl}=0$ in all other cases. Since $[\frac{\partial\tilde
x_k}{\partial y_l}]_{2n\times2n}=\varphi^{2-\frac{\lambda}{2}}(x_0)[O(1)]$ on $W(x_0)$ when $x_0$ near $E$,
\begin{equation}\label{eqn461}
\begin{split}
\Phi^{2-\frac{\lambda}{2}}(y) &=\varphi^{2-\frac{\lambda}{2}}(x_0)+(2-\frac\lambda2)
\Phi^{1-\frac{\lambda}{2}}\sum_{k,l=1}^{2n}D_{\tilde x_k}\varphi(\theta)
\frac{\partial\tilde x_k}{\partial y_l}y_l\\
&=\varphi^{2-\frac{\lambda}{2}}(x_0)+\varphi^{2-\frac{\lambda}{2}}(x_0)o(1)\\
&=\varphi^{2-\frac{\lambda}{2}}(x_0)(1+o(1))\\
&=\varphi^{2-\frac{\lambda}{2}}(x_0)O(1).
\end{split}
\end{equation}
where $\theta\in W(x_0)$ and $\Phi^{2-\frac{\lambda}{2}}(y)$
uniformly approximate the same quantity of
$\varphi^{2-\frac{\lambda}{2}}(x_0)$. Moreover,
\begin{equation}
\begin{split}
& D_{y_l}\left(\frac{\Phi^{2-\frac{\lambda}{2}}(y)}
{\varphi^{2-\frac{\lambda}{2}}(x_0)}\right)
=(2-\frac\lambda2)\frac{\Phi^{1-\frac{\lambda}{2}}(y)}
{\varphi^{2-\frac{\lambda}{2}}(x_0)}\sum_{k=1}^{2n}(D_{\tilde
x_k}\varphi)\left(\frac{\partial\tilde x_k}{\partial y_l}\right)
=\varphi^{1-\frac{\lambda}{2}}(x_0)O(1),\\
& D_{y_l}\left( \frac{\varphi^{\lambda}(x_0)}{\Phi^\lambda(y)}\right)
=(-\lambda)\varphi^{\lambda}(x_0)\Phi^{-\lambda-1}(y)\sum_{k=1}^{2n}(D_{\tilde x_k}\varphi)\left(\frac{\partial\tilde x_k}{\partial y_l}\right)
=\varphi^{1-\frac{\lambda}{2}}(x_0)O(1).
\end{split}
\end{equation}

Since $e_{k,2n-1}(x_0)=0$ when $k<n$, we have
\begin{equation}
\frac{E_{k,2n-1}(y)}{\varphi^{2-\frac{\lambda}{2}}(x_0)}
=\frac{e_{k,2n-1}(\tilde
x)-e_{k,2n-1}(0)}{\varphi^{2-\frac{\lambda}{2}}(x_0)}
=O(1).
\end{equation}
Thus we can write $e_{k,2n-1}(\tilde{x})=l_k^\prime(\tilde x)+
O(|\tilde{x}|^2),$ where $l_k^\prime$ is a linear function of $\tilde x$. It follows that
\begin{equation}
\lim\limits_{x_0\rightarrow E}
\left|\varphi^{-2+\frac{\lambda}{2}}(x_0)E_{k,2n-1}-l_k^\prime
(y_1,\ldots,y_{2n-2},0,0)\right|_{m,W(x_0)}=0.
\end{equation}
Similarly, by a direct calculation,
\begin{equation}
\frac{R_k(y)}{\varphi^{2-\frac{\lambda}{2}}(x_0)}\rightarrow 0
\end{equation}
when $x_0\rightarrow E, k=1,\ldots,n-1.$

Combining all the facts above, we conclude that if $k<n$
\begin{equation}\label{eqn435}
\begin{split}
&\lim\limits_{x_0\rightarrow E}\left|L_k-\left(\frac{\partial}{\partial y_{2k-1}}-\sqrt{-1}\frac{\partial}{\partial y_{2k}}+l_k\frac{\partial}{\partial y_{2n-1}}\right)\right|_{m,W(x_0)}=0,\\
\end{split}
\end{equation}
where $l_k=l_k^\prime(y_1,\ldots,y_{2n-2},0,0)$, and that
\begin{equation}\label{eqn436}
\begin{split}
\lim\limits_{x_0\rightarrow
E}\left|L_n-\left(\sqrt{-1}\frac{\partial} {\partial
y_{2n-1}}+\frac{\partial}{\partial
y_{2n}}\right)\right|_{m,W(x_0)}=0.
\end{split}
\end{equation}
Write $D=[d_{ij}]_{n\times 2n}$, $B=[b_{ij}]_{n\times 2n}$. Here
$d_{ij}$ and $b_{ij}$ are as defined in (ii) of the proposition.
Define
 $\tilde{D}=\left[\begin{array}{c}
              D \\
              \overline{D}
            \end{array}\right]_{2n\times 2n}
 $,
$\tilde{B}=[\begin{array}{cc}
             B^t & \overline{B}^t
           \end{array}]
$. Then $\tilde{D}\cdot\tilde{B}=I_{2n\times2n}$.
In order to prove $\{d_{ij}\}_{1\leq i\leq n, 1\leq j\leq 2n}$ and
the derivative of $\{d_{ij}\}_{1\leq i\leq n, 1\leq j\leq 2n}$ are
uniformly bounded, we only need to prove that the absolute value
$|det\wt{B}|$ of determinant of matrix $\tilde{B}$ has a uniform
lower bound. Let $A=[a_{ij}]_{2n\times2n}$. Here, $a_{2i-1,i}=1,
a_{2i,i}=-\sqrt{-1}, 1\leq i\leq n-1, a_{2n-1,n}=\sqrt{-1},
a_{2n,n}=1, a_{k,n+l}=\overline{a_{k,l}}$ and in other cases
$a_{ij}=0$. Let $C=[1]_{2n\times2n}$, $E=[e_{kl}]_{2n\times2n}$,
$e_{2n-1,k}=1, e_{2n-1,k+n}=1$ and  $e_{kl}=0$ in other cases. Then from
(\ref{eqn435}) and (\ref{eqn436}),
$\tilde{B}=A+O(\varphi^{2-\frac{\lambda}{2}}(x_0))C+O(\sigma_0)E$.
Thus $det(\tilde B)=detA\cdot
det(I+O(\varphi^{2-\frac{\lambda}{2}}(x_0))A^{-1}C+O(\sigma_0)A^{-1}E)$.
Since $|detA|=2^n\neq0$,  $|det\tilde{B}|$ has a uniform lower bound
with respect to $y\in W(x_0)$, when $\varphi(x_0)$ and $\sigma_0$
are sufficiently small. This proves $(ii)$.

Since $\{L_i\}_{i=1}^{n}$ is
an orthornormal basis with respect to the scaled Hermitian metric we have defined
on $W(x_0)$ near $E$, and from (\ref{eqn435}), (\ref{eqn436})
we see that the metric tensor and any order of its covariant differentiation on $W(x_0)\cap M$
induced from the Hermitian metric on $\-{S_{\varphi}^{\e_0}}\sm E$ near $E$
must have  uniform bounds. It is easy to see that there is a
constant $0<k_0<<1$ such that any ball in the induced metric over
$M_0\setminus E$ of radius $k_0\sigma_0'$ is contained in some $W(x_0)\cap
M$. From (\ref{eqn435}), (\ref{eqn436}) or by what we argued above,
the volume form in $W(x_0)\cap M$ in terms of the $y'$-coordinates
is uniformly bounded from above and from below by a positive
constant independent of $x_0$. We thus see that the volume of the
ball $B_{\mu \sigma'_0}$ in $M$ with radius $\mu \sigma'_0$ and
$\mu\le k_0$ has the estimate: $C_1(\mu\sigma^\prime_0)^{2n-1}\le
\mathrm{Vol}(B_{\mu \sigma'_0})\le C_2(\mu\sigma'_0)^{2n-1}$ with $C_1$ and
$C_2$ two positive constants.

Write $W'(x_0):=\{y': (y',0)\in W(x_0), |y'|<\sigma_0/2\}$. When
$\sigma_0$ is smaller than a certain fixed number, we can assume
that $W(x_0)\cap M$ is contained in $B_{\mu'\sigma_0'}$  and
$W'(x_0)$ contains $B_{\mu^{''}\sigma_0'}$, where
$0<\mu^{''}<\mu'<k_0$ are constants depending only on $\sigma_0$, but
not $x_0$. We can choose a family of $\{W(x_\alpha)\}$ near $E$ such
that (1). $\{W'(x_\alpha)\}$ covers $M_0\setminus E$ near $E$, (2). The
distance between any two centers in the induced metric over $M_0\setminus E$
is at least $\mu''\sigma_0'$. (The existence of such a family
follows from a simple construction based on the Zorn lemma). By the
just mentioned volume estimates, one concludes that such a cover is a
Besicovitch covering. Namely, there is a constant $N_0$, independent
of $\sigma_0$, such that any point is contained in at most
$N_0$-charts. Let $\wt{\xi_\a}\in C^{\infty}_0(W(x_\a)\cap M)$ be such
that $\wt{\xi_a}\equiv 1$ over $W'(x_\a)$ and the differentiation in
the $y'$-coordinates is bounded by $4/\sigma_0$. Define
$\xi_\alpha=\frac{\widetilde{\xi_a}}{\sqrt{\sum_{\alpha}\wt{\xi^2_a}}}.$ Then $\xi_\alpha\in
C^{\infty}_0(W(x_\alpha)\cap M)$, $\sum_{\alpha}{\xi^2_a}\equiv 1$ and
$|D_{y'}\xi_a|=O(\sigma^{-1}_0)$. This proves $(iv)$.

Finally, we note that if $L_j^v,\  j=1,\ldots,n-1, $ are replaced by
$$X_j^v=\sum_{k=1}^{n-1}U_{jk}L_k^v,$$ where $[U_{jk}]$ is a suitably chosen unitary
matrix such that
\begin{equation}
\partial\overline\partial t(X_i^v,\overline{X_j^v})|_{x_0}=d_{ij},1\leq i,j\leq n-1,
\end{equation}
where $d_{ij}=0,i\neq j; d_{jj}=d_j, j<n$.

Since $r_k=\frac{X_k^v(r)}{X_n(r)}=t\frac{X_k^v\varphi^{-4}(x)}{X_n^v(r)}$, we have
$r_k=0$ along $M\cap W(x_0)$. This shows that along $M\cap
W(x_0)$,
\begin{equation}\label{eqn25}
L_j=\varphi^{2-\frac{\lambda}{2}}X_j^v, j<n;~
\omega^n=\frac12\varphi^{-4+\lambda}(x)(dt+\sqrt{-1}\eta)
\end{equation}
Then
\begin{equation}\label{eqn462}
\begin{split}
\omega^n([L_k,\overline{L_n}])(x_0)
&=-d\omega^n(L_k,\overline{L_n})(x_0)\\
&=\varphi^{6-\frac32\lambda}(x_0)(-d\omega^n)(X_i^v, \overline{Z_n})(x_0)\\
&=-\varphi^{6-\frac32\lambda}(x_0)\left\{X_i^v(\omega^n(\overline{Z_n}))
-\overline{Z_n}\omega^n(X_i^v)-\omega^n([X_i^v, \overline{Z_n}])\right\}(x_0)\\
&=O(\varphi^{1-\frac{\lambda}{2}}(x_0))
\end{split}
\end{equation}
and
\begin{equation}
\begin{split}
\omega^n([L_i, \overline L_j])(x_0)
&=-d\omega^n(L_i, \overline L_j)(x_0)\\
&=-d\omega^n\left(\varphi^{2-\frac{\lambda}{2}}
(x_0)X^v_i(x_0),\varphi^{2-\frac{\lambda}{2}}(x_0)
\overline {X^v_j}(x_0)\right)\\
&=\varphi^{4-\lambda}(x_0)\omega^n([X_i^v, \overline{X_j^v}])(x_0)\\
&=\frac12(dt+\sqrt{-1}\eta)([X_i^v,\overline{X_j^v}])(x_0)\\
&=\partial t([X_i^v, \overline{X_j^v}])(x_0)\\
&=\partial\overline\partial t (X_i^v,\overline{X_j^v})(x_0).
\end{split}
\end{equation}
It follows that
\begin{equation}\label{26}
\omega^n([L_i,\overline{L_j}])(x_0)=d_{ij},1\leq i,j\leq n-1,
\end{equation}
and thus
\begin{equation}\label{eqn26}
\begin{split}
&\omega^n([L_i,\overline{L_j}])(y)=O(\sigma_0),~
i\neq j,1\leq i,j\leq n-1;\\
&\omega^n([L_j,\overline{L_j}])(y)=d_j+O(\sigma_0),~j<n
\end{split}
\end{equation}
in $W(x_0).$ Thus we obtain $(iii)$. The proof of the Proposition \ref{thm1} is
complete.
\end{proof}

Let $dV$ denote the volume form associated with the Hermitian metric
 defined before. In the coordinates $(y_1,\ldots,y_{2n})$ over
$W(x_0)$,  write $dV=V(y)dy$, where $dy=dy_1\cdots dy_{2n}$. Then we
have, as mentioned before, that  $V(y)$ satisfies the following:
\begin{equation}\label{eqn27}
|V(y)|_{1,W(x_0)}\leq a_1, \inf\limits_{y\in W(x_0)}V(y)>a_2>0,
\end{equation}
where $a_1, a_2$ are constants independent of $x_0$.

We will define inner product of two functions $g,h\in
C_c^{\infty}(\overline {S_\varphi^{\e_0}}\setminus E)$ by $$(g,h)=\int g\overline h dV.$$ Let
N be a submanifold of dimension $2n-1$ in $W(x_0)$ and let $ds$ be
the volume form of N that comes from Euclidean metric in
$(y_1,\ldots, y_{2n})$-variables. The following is the  divergence
theorem:

\begin{lemma}\label{lemma2}$(\text{Divergence theorem})$
Let $D\subset\subset \mathbb{R}^N$ be a smoothly bounded domain, $g,h,v\in C^{\infty}(\overline D)$, $V\neq 0$ on $\overline D, dV=Vdx, dS=Vds$. If $L=\sum\limits_{j=1}^nb_j\frac{\partial}{\partial x_j}$ is a smooth vector field,then
\begin{equation}\label{eqn28}
\int_D(Lg)\overline hdV=-\int_Dg\overline{\overline Lh}dV-\int_D eg\overline hdV+\int_{\partial D}g\overline h<L,{\bf{n}}>dS,
\end{equation}
where $e=\frac{LV}{V}+\sum\limits_{j=1}^N\frac{\partial
b_j}{\partial x_j}$, $\bf{n}$ is the unit outward   normal vector at
the boundary points, and $<,>$ is the Euclidean inner product in
$\mathbb{R}^N$.
\end{lemma}

Applying Lemma \ref{lemma2} 
one obtains the following: (See  [Lemma 5.7, Cat])
\begin{lemma}\label{lem3}
Let $L_1,\ldots,L_n$ be the frame  over $W(x_0)$ constructed in Proposition \ref{thm1}, then there exists functions $e_j\in C^{\infty}(W(x_0))$ and a function $P\in C^{\infty}(W(x_0))$
such that for all $g,h\in C_c^{\infty}(W(x_0))$
\begin{equation}\label{eqn29}
\begin{split}
&(L_jg,h)=-(g,\overline{L_j}h)-(e_jg,h),j=1,\ldots,n-1.\\
&(L_ng,h)=-(g,\overline{L_n}h)-(e_ng,h)+\int_{M_0}Pg\overline hdS-\int_{M_1}Pg\overline hdS,
\end{split}
\end{equation}
where $dS=Vds, M_0=\{z:r(z)=0\}, M_1=\{z:r(z)=-1\}.$ The real part
and imaginary part of the function $P$ satisfy: $0<c<Re(P(y))<C,$
$|Im(P(y))|\ll1$  for $y\in W(x_0)$ with $0<\varphi(x_0)\ll 1.$ Here
$c$ and $C$ are constants independent of $x_0$. Moreover,
$|e_j|_{m,W(x_0)}\leq C_m$ for $0<\varphi(x_0)\ll 1.$ Here $C_m$ is
a constant independent of $x_0$.
\end{lemma}
\begin{proof} Applying Lemma \ref{lemma2}, one
can see the above mentioned expressions hold with
\begin{equation}\label{eqn73}
e_j=\sum\limits_{k=1}^{2n}\frac{\partial b_{jk}}{\partial
y_k}+\sum\limits_{k=1}^{2n}b_{jk}\frac{\partial V(y)}{\partial
y_k}\frac{1}{V(y)},
\end{equation}
and
 $$P=L_n(y_{2n}) =\frac{\varphi^{\lambda}(x_0)}{\Phi^\lambda(y)}
-4\sum_{l=1}^{2n-1}y_{2n}\Phi^{3-\lambda}\Phi_lE_{nl}.$$ Since
$\frac{\varphi^{\lambda}(x_0)}{\Phi^\lambda(y)}$ approaches
uniformly to $1$ and
$\sum_{l=1}^{2n-1}y_{2n}\Phi^{3-\lambda}\Phi_lE_{nl}$ approaches
uniformly to zero as $x_0$ approaches to $E$, we conclude the proof
of the lemma.
\end{proof}

Now, suppose that when $0<\varphi(x_0)<\varepsilon_1<<1$, we have
constructed the special coordinates and special frame on $W(x_0)$ as
in Proposition \ref{thm1}. Notice that the subset
$${\cal K}_{\varepsilon_1}:=\-{S_\varphi^{\e_0}}\setminus
\left(\left\{x|\varphi(x)<\frac{\varepsilon_1}{2}\right\}\times
[-1,0]\right)$$ is compact in $S_\varphi^{\e_0}$. Here, as we
mentioned in $\S 2$, we identify ${\cal N}_{\e_0}$ with $\O_{\e_0}$.
We can then cover ${\cal K}_{\varepsilon_1}$ with finitely many
coordinate charts.

Write $\mathcal{E}^{0,q}(S_\varphi^{\e_0})$ for the set of smooth
$(0,q)$-forms  over  $S_\varphi^{\e_0}$. And write
$\mathcal{E}_c^{0,q}(S_\varphi^{\e_0})$ for the set of smooth
$(0,q)$-forms $U$ over $S_\varphi^{\e_0}$ with compact support in
$\-{S_\varphi^{\e_0}}\sm E$. Let
$\mathcal{E}_0^{0,q}(S_\varphi^{\e_0})$ denote the set of smooth sections
of $\Gamma^{(0,q)}(S_\varphi^{\e_0})$ with compact support in the
interior of $S_\varphi^{\e_0}$. We define the space
$L^2_{(0,q)}(S_\varphi^{\e_0})$ of $L^2$-integral $(0,q)$-forms by
using the scaled Hermitian metric.

Suppose $U\in\mathcal{E}_c^{(0,q)}(S_\varphi^{\e_0}),
U=\sum\limits_{|J|=q}U_J\overline{\omega^J}$. Then
\begin{equation}
\begin{split}
&\overline\partial
U=\sum\limits_{j=1}^n\sum\limits_{|J|=q}(\overline{L_j}U^J)
\overline\omega^j
\wedge\overline\omega^J+\cdots,\\
&\overline\partial^\prime
U=-\sum_{j=1}^n\sum\limits_{|K|=q-1}L_jU^{jK}
\overline{\omega^K}+\cdots,
\end{split}
\end{equation}
where $\overline\partial^{\prime}$ is the formal adjoint operator of
$\overline\partial$  and dots indicate terms where no derivatives of
$U^J$ occur. We extend $\overline\partial$ and
$\overline\partial^\prime$ to the $L^2$-space as in the
introduction. Then as in [Cat], define $\mathcal{B}^q$ for a
subspace of $\mathcal{E}_c^{(0,q)}(S_\varphi^{\e_0})$, whose
elements satisfy the $\overline\partial$-Dirichlet condition defined
in $\S 2$ along $M_0$ and the $\-{\p}$-Neumann condition along
$M_1$. Then $\mathrm{Dom}(T)\cap \mathrm{Dom}(S^\ast)\cap
\mathcal{E}_c^{(0,q)}(S_\varphi^{\e_0})=\mathcal{B}^q$. Moreover, as
in the [Cat] ( Lemma 6.4 of [Cat]), the H\"ormander-Friderichs
smooth lemma also holds in this setting:

Let $U\in \mathrm{Dom}(S^\ast)\cap \mathrm{Dom}(T)$. Then there exists $U_\mu\in
\mathcal{B}^q$ such that
\begin{equation}
\lim\limits_{\mu\rightarrow\infty}{\|U_\mu-U\|+\|S^\ast U_\mu-S^\ast
U\|+\|T U_\mu-TU\|}=0.
\end{equation}
Hence, in what follows, we need only to prove the estimate in our
main theorem for $U\in \mathcal{B}^q.$

\section{The $L^2$-estimate for the operator $T$ near the corner $E$}
In this section, we establish the estimate near $E$ for forms in
$\mathcal{B}^q.$ We follow the known
procedure to compute the the $Q$ norms as in [Ho1], [FK] and [Cat]. In
particular, we follow the computation in [Cat] and make
the needed modification to fit our situation here.

We first suppose $U\in \mathcal{B}^q$ with $\mathrm{supp}~U$ a
compact subset of $W(x_0)$  for some $x_0\in M_0$ with
$0<\varphi(x_0)\ll 1.$ Then
\begin{equation}\label{eqn33}
\begin{split}
&TU=\overline\partial
U=\sum\limits_{|J|=q}\sum\limits_{j=1}^n(\overline{L_j}U^J)
\overline\omega^{jJ}+\cdots,\\
&S^\ast
U=-\sum\limits_{|K|=q-1}\sum_{j=1}^{n}L_jU^{jK}\overline\omega^K+\cdots.
\end{split}
\end{equation}
Define $$Q(U,U)=\|TU\|^2+\|S^\ast U\|^2$$ and let
\begin{equation}\label{eqn34}
\begin{split}
&AU=\sum\limits_{|J|=q}\sum\limits_{j=1}^n(\overline{L_j}U^J)
\overline\omega^{jJ},\\
&BU=-\sum\limits_{|K|=q-1}\sum_{j=1}^{n}L_jU^{jK}\overline\omega^K.
\end{split}
\end{equation}
Immediately, we have
\begin{equation}\label{eqn35}
2\|S^\ast U\|^2+2\|TU\|^2+C_0\|U\|^2\geq\|AU\|^2+\|BU\|^2,
\end{equation}
where $C_0$ is a constant only depend on the coefficients of $L_j$
and $\omega_j$ and independent of $x_0$ and $U$.

Notice that
$$\|AU\|^2=\sum\limits_{(j, J)\neq (n, Kn)}^n\|\overline{L_j}U^J\|^2
-\sum\limits_{|K|=q-1}\sum\limits_{(j,k)\neq
(n,n)}(\overline{L_j}U^{kK},\overline{L_k}U^{jK})$$ where the
property that $(j,J)\neq(n,Kn)$  means that we  exclude those terms
where $j=n$ or $n\in J$.

We also notice that
\begin{equation}\label{eqn37}
\|BU\|^2=\sum\limits_{j,k=1}^n(L_jU^{jK}, L_kU^{kK}).
\end{equation}
To compute $\|AU\|^2+\|BU\|^2$, we follow the computation of Catlin
in [pp504-506, Cat] as follows:
First we calculate $\|BU\|^2$ by  several steps. Notice that if
$(j,k)\neq(n,n)$, $|K|=q-1$, then
\begin{equation}\label{eqn38}
\begin{split}
(L_kU^{kK},L_jU^{jK})&=(\overline{L_j}U^{kK},\overline{L_k}U^{jK})
+(\overline{e_j}L_kU^{kK},U^{jK})\\
&-(e_k\overline{L_j}U^{kK},U^{jK})
+([L_k,\overline{L_j}]U^{kK},U^{jK}).
\end{split}
\end{equation}
Define
\begin{equation}\label{eqn39}
L(U)=\sum\limits_{|K|=q-1}\|L_nU^{nK}\|^2+\sum\limits_{|J|=q,
n\not\in J}\|\overline{L_n}U^J\|^2
+\sum_{j=1}^{n-1}\sum\limits_{|J|=q}(\|L_jU^J\|^2+\|\overline{L_j}U^J\|^2).
\end{equation}
From (\ref{eqn9}) and by the standard big-small constant argument,
it follows that
\begin{equation}\label{eqn40}
\left|(\overline{e_j}L_kU^{kK},U^{jK})\right|\leq
\frac{C_1}{K_0}L(U)+C_1K_0\|U\|^2.
\end{equation}
Since
$$\left|(e_k\overline{L_j}U^{kK},U^{jK})\right|
=|(\overline{e_j}e_kU^{kK},
U^{jK})-(e_kU^{kK},L_jU^{jK})-(\overline{L_j}(e_k)U^{kK},
U^{jK})|,$$ we have
\begin{equation}\label{eqn444}
\left|(e_k\overline{L_j}U^{kK},U^{jK})\right|\leq
\frac{C_1}{K_0}L(U)+C_1K_0\|U\|^2.
\end{equation}
Here  $K_0, C_{j}s$ are constant independent of the choices of $x_0$
and $U$, which may be different in different contexts. $K_0$ is
supposed to be sufficiently large. Notice that
\begin{equation}\label{eqn41}
([L_k,\overline{L_j}]U^{kK},U^{jK})=\sum\limits_{i=1}^n(C_{kj}^iL_iU^{kK},U^{jK})
+\sum\limits_{i=1}^n(d_{kj}^i\overline{L_i}U^{kK},U^{jK}),
\end{equation}
where
\begin{equation} \label{a3}
C_{kj}^i=\omega^i([L_k,\overline{L_j}]),~
d_{kj}^i=\overline{\omega^i}([L_k,\overline{L_j}]).
\end{equation}
(a) If $i<n$,
then
\begin{equation}\label{eqn42}
|(C_{kj}^iL_iU^{kK},U^{jK})|+|(d_{kj}^i\overline{L_i}U^{kK},U^{jK})|\leq
\frac{C_1}{K_0}L(U)+C_1K_0\|U\|^2.
\end{equation}
(b) If $i=n, k=n$, and $j\neq n$
\begin{equation}\label{433}
|(C_{nj}^nL_nU^{nK},U^{jK})|\leq \frac{C_1}{K_0}L(U)+C_1K_0\|U\|^2,
\end{equation}
Hence, the remaining terms to be estimated include the following
\begin{equation}\label{43}
|(d_{nj}^n\overline{L_n}U^{nK},U^{jK})|.
\end{equation}
(c) If $i=n, j=n$, and $k\neq n$, we have
\begin{equation}\label{eqn43}
|(d_{kn}^n\overline{L_n}U^{kK},U^{nK})|\leq
\frac{C_1}{K_0}L(U)+C_1K_0\|U\|^2.
\end{equation}
Thus it suffices to estimate:
\begin{equation}\label{eqn44}
|(C_{kn}^nL_nU^{kK},U^{nK})|.
\end{equation}
As in [Cat], we use the following result on the standard uniform
sub-elliptic estimate to handle it: (See [CS], for instance.)
\begin{lemma}\label{lem7}
Let $\Omega$ be a bounded open neighborhood of the origin in
$\mathbb{R}^n$. Let $X_i, 1\leq i\leq k, k\leq n$ be vector fields
with smooth real coefficients up to $\-{\O}$. Let $\mathcal{L}_1$
 be collection of the $X_i^s, 1\leq i\leq k$ and $\mathcal{L}_2$
 be the union of $\mathcal{L}_1$ and the vectors of the form
$[X,Y]$ with $X, Y\in \mathcal{L}_1$. If $\mathcal{L}_2$ span the
tangent space of $\Omega$, then there exists $C>0$ such that
\begin{equation}\label{eqn m2}
\|u\|_{\frac12}^2\leq
C\left(\sum\limits_{i=1}^k\|X_iu\|^2+\|u\|^2\right),u\in
C_{0}^{\infty}(\Omega).
\end{equation}
Here C only depends on the coefficients of the vector fields.
\end{lemma}

For any $f\in C_c^{\infty}(W(x_0))$, we define the tangential Fourier transform for $f$ in $W(x_0)$ by
\begin{equation}
\hat f(\xi, y_{2n})=\int_{\mathbb{R}^{2n-1}}e^{-i<y^\prime, \xi>}f(y^\prime, y_{2n})dy_1\cdots dy_{2n-1},
\end{equation}
where $\xi=(\xi_1,\cdots, \xi_{2n-1})$ and $<y^\prime, \xi>=y_1\xi_1+\cdots+y_{2n-1}\xi_{2n-1}$.
We define the tangential Sobolev norm $|||f|||_s$ by
\begin{equation}
|||f|||_{s}^2=\int_{-\varphi^\lambda(x_0)}^0
\int_{\mathbb{R}^{2n-1}}|\hat
f(\xi,y_{2n})|^2(1+|\xi|^2)^{s}d\xi_1\cdots
d\xi_{2n-1}dy_{2n}, s\in \mathbb{R}.
\end{equation}
For more discussions  of the tangential Fourier transform and
tangential Sobolev norms, we refer the reader to [FK] and [CS] . From
Lemma \ref{lem7} and the uniform estimate of the coefficients of the
vector fields $\{L_i\}_{i=1}^n$ , there exists a constant $C_2$,
which does not depend on $x_0$ such that for all $f\in
C_c^\infty(W(x_0))$,

\begin{equation}\label{eqn54}
|||f|||_{\frac12}^2\leq C_2\sum\limits_{k=1}^{n-1}
\left(\|L_kf\|^2+\|\overline{L_k}f\|^2\right)+C_2^\prime\|f\|^2.
\end{equation}
\begin{lemma}\label{lem5}
Suppose $U\in \mathcal{B}^q(S_\varphi^{\e_0})$ and
$\text{supp~U}\subset\subset W(x_0)$ with $0<\varphi(x_0)\ll1.$ Then
\begin{equation}\label{eqn45}
|(C_{kn}^nL_nU^{kK},U^{nK})|\leq \frac{C_3}{K_0}L(U)+C_3K_0\|U\|^2,\
k\neq n,
\end{equation}
\begin{equation}\label{eqn455}
|(d_{nj}^n\overline{L_n}U^{nK},U^{jK})|\leq
\frac{C_3}{K_0}L(U)+C_3K_0\|U\|^2,\ j\neq n.
\end{equation}
\end{lemma}

\begin{proof}
We follow  the proof of [Lemma 7.8, Cat].
In [Cat], a useful fact is that the domain can be uniformly shrunk toward $M$, that helps
to get such types of estimates.
In our situation, the domain is fixed. However we go close and close
to $E$ such that the quantity
\begin{equation}\label{a2}
|C^n_{kn}|_{W(x_0)}=O(\varphi^{1-\frac{\lambda}{2}}(x_0)+\sigma_0)
\end{equation}
is sufficiently small when $0<\varphi(x_0)\ll1, \sigma_0\ll1$ to get the desired
estimates.
Write $L_n=(L_ny_{2n})(\overline{L_n}y_{2n})\overline{L_n}+\tilde
L_n$, then $\tilde L_n(y_{2n})=0$ and
\begin{equation}\label{eqn456}
\tilde L_n=\sum_{j=1}^{2n-1}\left\{L_n(y_j)-
\left(\frac{L_ny_{2n}}{\overline{L_n}y_{2n}}\right)
\overline{L_n}(y_j)\right\}\frac{\partial}{\partial y_j} ,
\end{equation}
\begin{equation}\label{eqn457}
\frac{L_ny_{2n}}{\overline{L_n}y_{2n}}=
\frac{1-4\sum_{l=1}^{2n-1}\frac{\Phi^3}{\varphi^{\lambda}(x_0)}y_{2n}E_{nl}(y)}
{1-4\sum_{l=1}^{2n-1}\frac{\Phi^3}{\varphi^{\lambda}
(x_0)}y_{2n}\overline{E_{nl}}(y)}.
\end{equation}
Let $\tilde L_n=\sum_{j=1}^{2n-1}a_j(y)\frac{\partial}{\partial
y_j}$. Then (\ref{eqn21}), (\ref{eqn462}), (\ref{eqn456}),
(\ref{eqn457}) give  that
\begin{equation}\label{eqn458}
\begin{split}
&|a_j(y)|_{n+1, W(x_0)}=O(\varphi^{2-\frac{\lambda}{2}}(x_0)), j<2n-1\\
&|a_{2n-1}(y)|_{n+1, W(x_0)}=O(1),
|C_{kn}^n|_{W(x_0)}=O(\varphi^{1-\frac{\lambda}{2}}(x_0)+\sigma_0).
\end{split}
\end{equation}
Then
\begin{equation}\label{eqn459}
(C_{kn}^nL_nU^{kK},
U^{nK})=\left(C_{kn}^n\frac{L_n(y_{2n})}{\overline{L_n}(y_{2n})}
\overline{L_n}U^{kK}, U^{nK}\right)+(C_{kn}^n\tilde L_nU^{kK},
U^{nK}).
\end{equation}
Now
\begin{equation}\label{eqn460}
\begin{split}
(C_{kn}^n\tilde L_nU^{kK}, U^{nK})&\leq
|C_{kn}^n|_{W(x_0)}\cdot|||\sum_{j=1}^{2n-1}a_j(y)
\frac{\partial U^{kK}}{\partial y_j}|||_{-\frac{1}{2}}\cdot|||U|||_{\frac{1}{2}}\\
&\leq C_3|C_{kn}^n|_{W(x_0)}\cdot|||U|||^2_{\frac12},
\end{split}
\end{equation}
On the other hand, since
$\|\overline{L_n}U^{kK}\|\leq L(U)$ then
$$\left(C_{kn}^n\frac{L_n(y_{2n})}{\overline{L_n}(y_{2n})}
\overline{L_n}U^{kK},
U^{nK}\right)\leq\frac{C_3}{K_0}L(U)+C_3K_0\|U\|^2, k<n.$$
Notice that we can make $|C_{kn}^n|_{W(x_0)}$ sufficiently small by
letting $x_0$ close to $E$.
Combining (\ref{eqn54}), (\ref{eqn460}) and the just obtained estimate, we conclude the
estimate  in (\ref{eqn45}). The proof for (\ref{eqn455}) is similar.
\end{proof}
(d) If $i=n$, and $1\leq j,k\leq n-1$, we need to control the
following terms:
$$(C_{kj}^nL_nU^{kK},U^{jK}),(d_{kj}^n\overline{L_n}U^{kK},U^{jK}) .$$
For $n\in K$, it holds that
\begin{equation}\label{eqn46}
|C_{kj}^n L_nU^{kK},U^{jK})|\leq \frac{C_2}{K_0}L(U)+C_2K_0\|U\|^2.
\end{equation}
When $n\not\in K$, we have
\begin{equation}\label{eqn445}
|(d^n_{kj}\overline{L_n}U^{kK},U^{jK})|\leq\frac{C_2}{K_0}L(U)+C_2K_0\|U\|^2.
\end{equation}
The only remaining two cases are  (i): For $n\not\in K$, we need to
control
$$(C_{kj}^nL_nU^{kK},U^{jK})$$ and (ii): For $n\in K$, we need to
control
$$(d^n_{kj}\overline{L_n}U^{kK},U^{jK}).$$
Since $C^n_{kj}(y)=\omega^n([L_k,\overline{L_j}])(y)
=d_{kj}(x_0)+O(\sigma_0),  d^n_{kj}(y)=-d_{kj}(x_0)+O(\sigma_0)$,
thus
\begin{equation}\label{eqn47}
(C_{kj}^nL_nU^{kK},U^{jK})=d_{kj}(x_0)(L_nU^{kK},U^{jK})
+(O(\sigma_0)L_nU^{kK},U^{jK}).
\end{equation}
Define $$E(x_0,U,
K_0)=\frac{1}{K_0}L(U)+K_0\|U\|^2+\sigma_0\int_{M_0}|U|^2 dS
+\sigma_0\int_{M_1}|U|^2 dS.$$
When $n\not\in K$, integrating by part, we have from (\ref{eqn47})
\begin{equation}\label{eqn48}
\begin{split}
&(C^n_{jj}L_nU^{jK},U^{jK})=-d_j(x_0)
\int_{M_1}P|U^{jK}|^2dS+O(E(x_0,U,K_0)),\\
&(C^n_{kj}L_nU^{kK},U^{jK})=O(E(x_0,U,K_0)).
\end{split}
\end{equation}
When $n\in K$, integrating by part, we have from (\ref{eqn47})
\begin{equation}\label{eqn565}
\begin{split}
&(d^n_{jj}\overline{L_n}U^{jK},U^{jK})=-d_j(x_0)\int_{M_0}P|U^{jK}|^2dS+O(E(x_0,U,K_0)),\\
&(d^n_{kj}L_nU^{kK},U^{jK})=O(E(x_0,U,K_0)).
\end{split}
\end{equation}
Then \begin{equation}\label{eqn555}
\begin{split}
\|AU\|^2+\|BU\|^2 &=\sum\limits_{|J|=q,n\not\in
J}\|\overline{L_n}U^{J}\|^2+\sum\limits_{|K|=q-1}\|L_nU^{nK}\|^2+
\sum_{j=1}^{n-1}\sum\limits_{|J|=q}\|\overline{L_j}U^{J}\|^2\\
&-\sum\limits_{|J|=q}\sum_{j\in J}\left(d_j(x_0)
\int_{M_0}P|U^J|^2dS+d_j(x_0)\int_{M_1}P|U^J|^2dS\right)\\
&+O(E(x_0,U,K_0))\\
&=\sum\limits_{|J|=q,n\not\in
J}\|\overline{L_n}U^{J}\|^2+\sum\limits_{|K|=q-1}\|L_nU^{nK}\|^2+
\sum_{j=1}^{n-1}\sum\limits_{|J|=q}\|\overline{L_j}U^{J}\|^2\\
&-\sum\limits_{|J|=q}\sum_{j\in J}\left(d_j(x_0)
\int_{M_0}Re(P)|U^J|^2dS+d_j(x_0)\int_{M_1}Re(P)|U^J|^2dS\right)\\
&+O(E(x_0,U,K_0)).
\end{split}
\end{equation}
Making use of the sign condition on the Levi forms and by a standard
argument (see [pp 62, FK] and [CS], for instance), we obtain from the above
\begin{equation}\label{eqn51}
\| AU\|^2+\|BU\|^2\geq c\left(L(U)+\int_{M_0}|U|^2
dS+\int_{M_1}|U|^2 dS\right) +O(K_0\|U\|^2).
\end{equation}
Hence, there exist positive constants $c$ and $C$ both
indpendent of the choices of $x_0$ and $ U$ such that
\begin{equation}\label{eqn52}
\|AU\|^2+\|BU\|^2\geq c L(U)-C\|U\|^2.
\end{equation}

The following Lemma from [Cat] is a fundamental fact, by which the
mixed boundary conditions enters the estimate.  It is Lemma 7.7 of
[Cat] with $\sigma^3$ being replaced by $\varphi^{\lambda}(x_0)$.

\begin{lemma}\label{lem6}
Suppose $f\in C_c^\infty(W(x_0))$ with $\varphi(x_0)\ll1$, and $f$
vanishes either on $M_0$ or $M_1$. Then there exists a constant
$\hat C_1$ independent of $x_0$ and $U$, so that
\begin{equation}\label{eqn102}
\begin{split}
&\varphi^{-\lambda}(x_0)\|f\|^2\leq
\hat C_1\left(\|L_nf\|^2+\sum\limits_{k=1}^{n-1}
\left(\|L_kf\|^2+\|\overline{L_k}f\|^2\right)\right)\\
&\varphi^{-\lambda}(x_0)\|f\|^2\leq \hat C_1\left(\|\overline
{L_n}f\|^2+\sum\limits_{k=1}^{n-1}
\left(\|L_kf\|^2+\|\overline{L_k}f\|^2\right)\right)
\end{split}
\end{equation}
\end{lemma}

Combining this Lemma \ref{lem6} with (\ref{eqn52}), we  proved Part
(1) of the following theorem:
\begin{theorem}\label{thm2} (1).
There exists a constant $0<\varepsilon_2\ll1$ independent of $x_0$
and a constant $\tilde{C}$ independent of $x_0$ and $\varepsilon_2$
such that if ~$0<\varphi^\lambda(x_0)\leq\varepsilon_2$ and $U\in
\mathcal{B}^q$ with $supp~U\subset\subset W(x_0)$, then
\begin{equation}\label{eqn66}
\varphi^{-\lambda}(x_0)\|U\|^2\leq \tilde{C}\left(\|TU\|^2+\|S^\ast
U\|^2\right).
\end{equation}

\noindent  (2). There exists a small neighborhood $V_c$ of $E$ in
$\-{S_\varphi^{\e_0}}$ such that for any ~$U\in \mathcal{B}^q $ it
holds that
\begin{equation}\label{eqn68}
\int_{V_c}|U|^2dV\leq \|TU\|^2+\|S^\ast U\|^2+\int_{S_\varphi^{\e_0}\setminus V_c}|U|^2dV.
\end{equation}
\end{theorem}
\begin{proof}
Let $V_c=\{(x,t)\in \mathcal{N}_{\e_0}\cap
S^{\e_0}_\varphi:0<\varphi^\lambda(x)\leq\varepsilon_3\}$, $\varepsilon_3<\varepsilon_2$. Let
$\{W(x_\alpha)\}$, $\{\xi_\alpha\}$ be as in Proposition 3.1 $(iv)$,
where $\{W(x_\alpha)\}$ is   a Besicotvich cover of $V_c$ for
$\varepsilon_3$ sufficiently small.
\begin{equation}\label{eqn70}
\begin{split}
\int_{V_c}|U|^2dV&=\sum\limits_\alpha\int_{V_c}|\xi_\alpha U|^2dV\le
\sum\limits_\alpha \int_{W(x_\alpha)}|\xi_\alpha
U|^2\leq\tilde{C}\sum\limits_\alpha \varphi^\lambda(x_\alpha)
(\|T\xi_\a U\|^2+\|S^\ast \xi_\alpha U\|^2)\\
&\leq\tilde{C}\sum\limits_\a\varphi^\lambda(x_\a) (\|\xi_\alpha
TU\|^2+\|\xi_\alpha S^\ast
U\|^2)+\tilde{C}\frac{c^2}{\sigma_0^2}\sum\limits_\alpha
\varphi^{\lambda}(x_\alpha)
\int_{W(x_\alpha)}| U|^2dV\\
&\leq\tilde{C}\varepsilon_3\left(\|TU\|^2+\|S^\ast
U\|^2\right)+\tilde{C}\varepsilon_3\frac{c^2}{\sigma_0^2}
N_0\int_{S_\varphi^{\e_0}}|U|^2dV.
\end{split}
\end{equation}
Here the $N_0$ is as in Proposition 3.1 $(iv)$. When $\varepsilon_3$ is
sufficiently  small such that
$$\max\left\{\tilde{C}\varepsilon_3,\tilde{C}
\varepsilon_3\frac{c^2}{\sigma_0^2} N_0\right\} \leq\frac12,$$ then we get
\begin{equation}\label{eqn81}
\int_{V_c}|U|^2dV\leq \|TU\|^2+\|S^\ast
U\|^2+\int_{S_\varphi^{\e_0}\setminus V_c}|U|^2dV.
\end{equation}
\end{proof}

\section{Proof of Theorem \ref{101}}
We now give a proof of Theorem \ref{101}. We use the notations set
up above. By the H\"ormander-Friderichs approximation theorem
mentioned at the end of $\S 3$, we need only work on forms in ${\cal
B}^q$.

 Let
$V_{c_1}=\{(x,t)\in \mathcal{N}_{\varepsilon_0}\cap
S^{\e_0}_\varphi:\varphi^\lambda(x)< \varepsilon_4\}$ with
$\varepsilon_4<\varepsilon_3$ be a smaller neighborhood of the
corner $E$ in $\overline{S^{\e_0}_\varphi}$ with $V_{c_1}\subset\subset V_c$. Let
$$M_{c,0}=M_0\setminus V_c,
M_{c_1,0}=M_0\setminus V_{c_1}.$$ Suppose that there exists a
tubular neighborhood $\mathcal{N}_{c_1}$ of $M_{c_1,0}$ in $X$ and a
$C^\infty$ map $\Phi_{c_1}$ such that
$\Phi_{c_1}:\mathcal{N}_{c_1}\rightarrow M_{c_1,0}\times(-2,2)$ is a
diffemorphism. Let
$\mathcal{L}_{c_1}={\Phi_{c_1}}_{\ast}(T^{1,0}\mathcal{N}_{c_1})$,
where $T^{1,0}\mathcal{N}_{c_1}$ is the holomorphic tangent bundle
of $\mathcal{N}_{c_1}$. Write $\Omega_{c_1}=M_{c_1}\times(-2,2)$,
then $(\Omega_{c_1},\mathcal{L}_{c_1})$ is a complex manifold
biholomorphic to $(\mathcal{N}_{c_1}, T^{1,0}\mathcal{N}_{c_1})$.
Also $M_{c_1,0}$ is a hypersurface in
$(\Omega_{c_1},\mathcal{L}_{c_1})$. In what follows, as before, when
there is no risk of causing confusion, we identify
$\mathcal{N}_{c_1}$ with $\Omega_{c_1}$ and objects defined over
$\mathcal{N}_{c_1} $ with those corresponding to $\Omega_{c_1}$. We
define two subdomains of $\mathcal{N}_{c_1}$ as follows:

\begin{equation}\label{eqn74}
\begin{split}
&\tilde S_{c,\varepsilon}=M_{c,0}\times[-\varepsilon^3,0],
\tilde S_{c_1,\varepsilon}=M_{c_1,0}\times[-\varepsilon^3,0],\\
&M_{c,\varepsilon}:=M_{c,0}\times\{-\varepsilon^3\},
M_{c_1,\varepsilon}:=M_{c_1,0}\times\{-\varepsilon^3\}.
\end{split}
\end{equation}
We can assume that $\tilde S_{c_1,\varepsilon}$ is contained in
$\mathcal{N}_{c_1}\cap S^{\e_0}_\varphi$ by making $\varepsilon$
sufficiently small. First, we  choose a finite cover
$\{V_v\}_{v=1}^m$ of $M_{c_1,0}$. With the same argument as in
Proposition 3.1 $(iv)$ we can assume that for any $x_0\in M_{c,0}$, there is
a coordinate neighborhood
\begin{equation}
V(x_0)=\{(x^\prime, x_{2n})\in M_{c_1,0}\times[-\varepsilon^3,0]: |x_1|< \varepsilon, \cdots, |x_{2n-1}|<\varepsilon, -\varepsilon^3\leq x_{2n}\leq 0\}
\end{equation}
with $V(x_0)$ contained in a certain $V_v$
when $\varepsilon$ is sufficiently small, where $x_0$ corresponds to the origin and $(x_1,\cdots,x_{2n-1})$ is the coordinates of $V(x_0)\cap M_{c_1,0}$. There exist a special orthnormal frame $\{L_i\}_{i=1}^n$ and its dual frame
$\{\omega_i\}_{i=1}^n$ satisfying the type of estimates as in Proposition 3.1 (iii).

Let $U\in\mathcal{E}_c^{(0,q)}(S^{\e_0}_\varphi)$ with compact
support in one of the above constructed coordinate neighborhoods
$V(x_0)$ such that it satisfies the $\overline\partial$-Dirichlet
condition on $M_{c_1,0}$ and vanishes near $M_{c_1,\varepsilon}$.
Then we have the following estimate
\begin{lemma}\label{lem00}
There exists a constant $\tilde{C_1}$ independent of $x_0$ and $U$ such that
\begin{equation}\label{eqna1}
\varepsilon^{-3}\|U\|^2\leq \tilde{C_1}Q(U,U)
\end{equation}
for $\varepsilon$  sufficiently small.
\end{lemma}
In order to prove Lemma \ref{lem00}, we first recall the definition of condition $Z(q)$ from [FK] and [CS].
\begin{Definition}
Let $D$ be a relatively compact subset with $C^\infty$ boundary in a
complex Hermitian manifold of complex dimension $n\geq2$. We say the
boundary $\p D$ satisfies condition $Z(q)$ if for $1\leq q\leq n-1$,
the Levi form associated with $D$ has at least $n-q$ positive
eigenvalues or at least $q+1$ negative eigenvalues at every boundary
point.
\end{Definition}

{\bf Proof of Lemma \ref{lem00}}:
We will use the
 property that $U$ vanishes near $M_{c_1, \varepsilon}$ to obtain the estimate in (\ref{eqna1}).
  We define the Hodge star operator $\star$ with respect to the scaled metric  as follows: For any $U_1, U_2\in\mathcal{E}_c^{p,q}(S^{\epsilon_0}_\varphi)$, we have
\begin{equation}
(U_1, U_2)=\int_{S^{\epsilon_0}_\varphi}U_1\wedge\star U_2 dV.
\end{equation}
Since $U\in\mathcal{E}_c^{(0,q)}(S_\varphi^{\epsilon_0})$ with $\mathrm{supp} U\subset\subset V(x_0)$ and
 $U$ satisfies the $\overline\partial$-Dirichlet condition on $M_{c_1,0}$, vanishes near $M_{c_1,\varepsilon}$,
 $\star U$ will be a $(n,n-q)$-form which satisfies the $\overline\partial$-Neumann condition
  on $M_{c_1,0}$ and vanishes near $M_{c_1, \varepsilon}$. Let $V=\star U$. We then have $V=\sum\limits_{|J|=n-q}V^J\omega^1\wedge\cdots\wedge\omega^n\wedge\overline\omega^J$. Then
\begin{equation}\label{a4}
\begin{split}
\|\overline\partial V\|^2+\|\overline\partial^\ast V\|^2
=&\sum\limits_{j, |J|=n-q}\|\overline L_j V^J\|^2+\sum\limits_{j,k, L}\int_{M_0\cap V(x_0)}C^n_{jk}PV^{jL}\overline{V^{kL}}dS\\
&+ O(l(V)\cdot\|V\|)+O(\|V\|^2)
\end{split}
\end{equation}
where $l(V)=\sum\limits_{j, J}\|\overline L_j V^J\|$ and
$C^n_{jk}=w^n([L_j, \overline L_k])$. It follows that there exist
constant $a_1, a_2$ which are independent of $x_0$ and $V(x_0)$ such
that
\begin{equation}
\begin{split}
\|\overline\partial V\|^2+\|\overline\partial^\ast V\|^2
\geq a_1l(V)-a_2\|V\|^2
&+\sum\limits_{|J|=n-q}\sum\limits_{j\in J}\int_{\p M_{c_1,0}}d_j(x_0)P|V^J|^2dS\\
&+O(\varepsilon)\sum\limits_{|J|=n-q}\sum\limits_{j\in J}\int_{\p
M_{c_1,0}}d_j(x_0)P|V^J|^2dS,
\end{split}
\end{equation}
where $\{d_j(x_0)\}_{j=1}^{n-1}$ are the eigenvalues of the
Levi-form on $M_0$. By the assumption of $M_0$, we know that $M_0$
satisfies Condition $Z(n-q)$. Proceeding in the standard way as in
[FK] and [CS], we see that there exist constant $a_3, a_4$ such that
\begin{equation}\label{a8}
\|\overline\partial V\|^2+\|\overline\partial^\ast V\|^2\geq a_3\left(\sum\limits_{j,J}\|\overline L_j V^J\|^2+\sum\limits_{j=1}^{n-1}\sum\limits_{J}\|L_j V^J\|^2\right)-a_4\|V\|^2.
\end{equation}
Since $V$ vanish near $M_{c_1,\varepsilon}$,  from Lemma \ref{lem6} we have
\begin{equation}\label{a9}
\sum\limits_{j,J}\|\overline L_j V^J\|^2+\sum\limits_{j=1}^{n-1}\sum\limits_{J}\|L_j V^J\|^2\geq \varepsilon^{-3}\|V\|^2.
\end{equation}
Combing (\ref{a8}) and (\ref{a9}) and when $\varepsilon$ is sufficiently small, we have
\begin{equation}\label{a99}
\|\overline\partial V\|^2+\|\overline\partial^\ast V\|^2\geq a_5\varepsilon^{-3}\|V\|^2
\end{equation}
when $\varepsilon$ is sufficiently small.
Substituting $V=\star U$ to (\ref{a99}), we have
\begin{equation}
\|\overline\partial \star U\|^2+\|\overline\partial^\ast \star U\|^2\geq a_5\varepsilon^{-3}\|\star U\|^2
\end{equation}
and since the Hodge star operator $\star$ is an isometry operator in $L^2$-space, we have
\begin{equation}\label{a09}
\|\star\overline\partial \star U\|^2+\|\star\overline\partial^\ast \star U\|^2\geq a_5\varepsilon^{-3}\|\star U\|^2
\end{equation}
Substituting the identity
$\overline\partial^\ast=-\star\overline\partial\star$ and $
\overline\partial=\star\overline\partial^\ast\star$ to (\ref{a09}),
we complete the proof of the Lemma {\ref{lem00}}. (For detailed
discussions  of the Hodge star operator in the $L^2$-space, we also
refer the reader to [CSh].) $\endpf$

\medskip
Moreover, we can choose a set $I$ such that $\{V(x_i)\}_{i\in I}$ is
a covering of $\tilde{S}_{c,\varepsilon}$. Moreover,
$V(x_i)\subset\subset \tilde{S}_{c_1,\varepsilon}$ and there exists
an integer $\hat N$, independent of $\varepsilon$, such that no
point of $\tilde{S}_{c,\varepsilon}$ lies in more than $\hat N$ of
such $V(x_i)$'s. We choose functions $\rho_i\in C_c^\infty(V(x_i))$
such that $\sum_{i\in I}\rho_i^2=1$ in a neighborhood of
$\tilde{S}_{c,\varepsilon}$ and $|\rho_i|_{1,V(x_i)}\leq \tilde
c_1\varepsilon^{-1}$, where $\tilde{c}$ is a constant independent of
$x_i$. Then for all
$U\in\mathcal{E}_c^{(0,q)}(S^{\epsilon_0}_\varphi)$ which satisfies
the  $\overline\partial$-Dirichlet condition on $M_{c_1,0}$ and
vanishes near $M_{c_1,\varepsilon}$, we have
\begin{equation}\label{eqn76}
\begin{split}
\varepsilon^{-3}\int_{\tilde S_{c,\varepsilon}}|U|^2dV &\leq
\varepsilon^{-3}\sum\limits_{i\in I}\int_{V(x_i)}|\rho_i U|^2dV
\leq \tilde{C_1}\sum\limits_{i\in I} Q(\rho_iU,\rho_iU)\\
&\leq
\tilde{C_1}Q(U,U)+\frac{\tilde{C_2}\tilde c_1^2\hat{N}}{\varepsilon^2}
\int_{\tilde{S}_{c_1,\varepsilon}}|U|^2dV
\end{split}
\end{equation}
Thus
\begin{equation}\label{eqn77}
\int_{\tilde S_{c,\varepsilon}}|U|^2dV\leq \tilde{C_1}\varepsilon^3Q
(U,U)+\tilde{C_2}\tilde c_1^2\hat{N}\varepsilon\int_{\tilde{S}_{c_1,\varepsilon}}|U|^2dV
\end{equation}
and
\begin{equation}\label{eqn556}
(1-\tilde C_2\tilde c_1^2\hat{N}\varepsilon)\int_{\tilde
S_{c,\varepsilon}}|U|^2dV\leq \tilde C_1\varepsilon^3Q(U,U)+\tilde
C_2\tilde c_1^2\hat{N}\varepsilon\int_{\tilde
S_{c_1,\varepsilon}\setminus\tilde S_{c,\varepsilon}}|U|^2dV
\end{equation}
First, we choose $\varepsilon$ such that $1-\tilde
C_2\tilde c_1^2\hat{N}\varepsilon\geq\frac12$, then
\begin{equation}\label{eqn557}
\begin{split}
\int_{\tilde S_{c,\varepsilon}}|U|^2dV
&\leq 2\tilde
C_1\varepsilon^3Q(U,U)+2\tilde
C_2\tilde c_1^2\hat{N}\varepsilon\int_{\tilde S_{c_1,\varepsilon}\setminus
\tilde S_{c,\varepsilon}}|U|^2dV \\
&\leq2\tilde
C_1\varepsilon^3Q(U,U)+2\tilde
C_2\tilde c_1^2\hat{N}\varepsilon\int_{V_c}|U|^2dV
\end{split}
\end{equation}
Next, let $$M_{c,1}= M_1\cap(S^{\e_0}_{\varphi}\setminus V_c),
M_{c_1,1}= M_1\cap (S^{\e_0}_{\varphi}\setminus V_{c_1}).$$ By the
same construction as above we can find a small  neighborhood
$\mathcal{T}_{c_1,\varepsilon} $ of $M_{c_1,1}$ in
$S^{\e_0}_\varphi$ which is biholomorphic to
$M_{c_1,1}\times[0,\varepsilon^3]$. We similarly write
\begin{equation}\label{eqn558}
\begin{split}
&\mathcal{T}_{c,\varepsilon}=M_{c,1}\times[0,\varepsilon^3], \mathcal{T}_{c_1,\varepsilon}=M_{c_1,1}\times[0,\varepsilon^3]\\
&\widehat{M}_{c,0}=M_{c,1}\times\{0\}, \widehat{M}_{c,\varepsilon}=M_{c,1}\times\{\varepsilon^3\}\\
&\widehat{M}_{c_1,0}=M_{c_1,1}\times\{0\},
\widehat{M}_{c_1,\varepsilon}=M_{c_1,1}\times\{\varepsilon^3\}
\end{split}
\end{equation}
For all $x_0\in \widehat{M}_{c,0}$ we  choose a coordinate neighborhood
$V(x_0)=\{(x^\prime, x_{2n})\in M_{c_1,0}\times[0,\varepsilon^3]:
x^\prime=(x_1,\cdots, x_{2n-1}), |x_1|<\varepsilon, \cdots,
|x_{2n-1}|<\varepsilon, 0\leq x_{2n}\leq \varepsilon^3\}$, where
$x_0$ corresponds to the origin and $(x_1,\cdots,x_{2n-1})$ is the
coordinates of $V(x_0)\cap \widehat{M}_{c,0}$. There exist a special
frame $\{L_i\}_{i=1}^n$ and their dual frame $\{\omega_i\}_{i=1}^n$,
that satisfy good estimates as in Proposition 3.1 with
$L_n(x_{2n})=-1$ on $V(x_0)\cap M_1$. Moreover,
\begin{equation}\label{eqn566}
\omega^n([L_i,\overline{L_j}])(x)=d_{ij}(x_0)+O(\varepsilon), 1\leq,
i,j\leq n-1,
\end{equation}
where $d_{ij}(x_0)=0$ when $i\neq j$, $d_{ij}(x_0)=d_j(x_0)$ when
$i=j$ and $\{d_j(x_0)\}_{j=1}^{n-1}$ are the  Levi eigenvalues  of the Levi form on $M_1$ with
respect to the domain. Notice that we now have  $P=-L_n(x_{2n})=1$
on $V(x_0)\cap M_1$.

By assumption of $M_1$, we have that $M_1$ satisfies $Z(q)$ condition.
 Let $U\in\mathcal{E}_c^{0,q}(S^{\e_0}_\varphi)$  satisfy the  $\overline\partial$-Neumann
condition on $\widehat{M}_{c_1,0}$ and vanish near $\widehat{M}_{c_1,1}$. Assume that $U$ has a compact
support in some $V(x_0)$. Then by a similar argument as in Lemma \ref{lem00} and when $\varepsilon$ is sufficiently small we have
the following $L^2$-estimate
\begin{equation}\label{eqn570}
\varepsilon^{-3}\|U\|^2\leq\tilde{C_3}Q(U,U).
\end{equation}
Here the  constant $\tilde{C_3}$  is independent of $U$ and $x_0$.
Then choosing a covering $\{V(x_\lambda)\}_{\lambda\in \Lambda}$ of $\mathcal{T}_{c,\varepsilon}$ with the same  property as  the covering $\{V(x_i)\}_{i\in I}$ in the proof of Lemma \ref{lem00} and also using a partition of
unity with the same property with respect to such
covering, we have
\begin{equation}\label{eqn560}
\int_{\mathcal T_{c,\varepsilon}}|U|^2dV\leq 2\tilde
C_3\varepsilon^3Q(U,U)+2\tilde
C_4\tilde c_2^2\hat{N}\varepsilon\int_{\mathcal
T_{c_1,\varepsilon}\setminus\mathcal T_{c,\varepsilon}}|U|^2dV
\leq2\tilde C_3\varepsilon^3Q(U,U)+2\tilde
C_4\tilde c_2^2\hat{N}\varepsilon\int_{V_c}|U|^2dV
\end{equation}
for some constants $\tilde{C_3},\tilde{C_4}, \tilde c_2$ which do not depend
on $U$.

We choose a cut-off function $\chi_{\varepsilon}$ such that
$\chi_{\varepsilon}\equiv1$ in a small neighborhood $O_\varepsilon$ of
$M_{c_1,0}\cup{M}_{c_1,1}\cup V_{c_1}$ and $\chi_{\varepsilon}$ equals to
zero near $M_{c,\varepsilon}$ and $\widehat{M}_{c,\varepsilon}$. Then there is a neighborhood
$V_{c,0,\varepsilon}$ of $M_{c,0}$ in $\tilde{S}_{c,\varepsilon}$
such that when $\chi_{\varepsilon}$ restricted to
$V_{c,0,\varepsilon}$ equals to 1. Then
\begin{equation}\label{eqn78}
\begin{split}
\int_{V_{c,0,\varepsilon}}|U|^2dV
&\leq\int_{\tilde{S}_{c,\varepsilon}}|\chi_\varepsilon U|^2
\leq2\tilde{C_1}\varepsilon^3 Q(\chi_\varepsilon U,\chi_\varepsilon U)+2\tilde{C_2}\tilde c_1^2\hat{N}\varepsilon\int_{V_c}|\chi_\varepsilon U|^2dV\\
&\leq2\tilde{C}_1\varepsilon^3Q(U,U)+C(\varepsilon)\int_{K_\varepsilon}|U|^2dV
+2\tilde C_2\tilde c_1^2\hat{N}\varepsilon\int_{V_c}|U|^2dV,
\end{split}
\end{equation}
where $C(\varepsilon)$ is a constant depending on $\varepsilon$ and
$K_\varepsilon$ is the compliment of $O_\varepsilon$ in the domain, which is
a compact subset of $S^{\e_0}_\varphi$.

Similarly, there is a neighborhood $V_{c,1,\varepsilon}$ of
${M}_{c,1}$ in $\mathcal T_{c,\varepsilon}$ such that
$\chi_\varepsilon$ is identically one when restricted to
$V_{c,1,\varepsilon}$. Then
\begin{equation}\label{eqn778}
\begin{split}
\int_{ V_{c,1,\varepsilon}}|U|^2dV
&\leq\int_{\mathcal{T}_{c,\varepsilon}}|\chi_\varepsilon U|^2
\leq2\tilde{C_3}\varepsilon^3 Q(\chi_\varepsilon U,\chi_\varepsilon U)+2\tilde{C_4}\tilde c_2^2\hat{N}\varepsilon\int_{V_c}|\chi_\varepsilon U|^2dV\\
&\leq2\tilde{C}_3\varepsilon^3Q(U,U)+C(\varepsilon)\int_{K_\varepsilon}|U|^2dV
+2\tilde C_4\tilde c_2^2\hat{N}\varepsilon\int_{V_c}|U|^2dV.
\end{split}
\end{equation}
Thus
\begin{equation}\label{eqn83}
\begin{split}
&\int_{V_c}|U|^2dV+\int_{V_{c,0,\varepsilon}}|U|^2dV+\int_{ V_{c,1,\varepsilon}}|U|^2dV\\
&\leq(1+4\tilde{C}_1\varepsilon^3+4\tilde{C}_3\varepsilon^3)Q(U,U)
+4C(\varepsilon)\int_{K_\varepsilon}|U|^2dV\\
&+\int_{K_\varepsilon^\ast}|U|^2dV+(4\tilde
C_2\tilde c_1^2+4\tilde{C}_4\tilde c_2^2)\hat{N}\varepsilon\int_{V_c}|U|^2dV,
\end{split}
\end{equation}
where $K_\varepsilon^\ast=S^{\e_0}_\varphi\setminus\{V_c\cup
V_{c,0,\varepsilon}\cup V_{c,1,\varepsilon}\}$ which is a compact
subset of $S^{\e_0}_\varphi$. We choose $\varepsilon$ such that
$(4\tilde C_2\tilde c_1^2+4\tilde{C}_4\tilde
c_2^2)\hat{N}\varepsilon\leq\frac12$. Then
\begin{equation}\label{eqn83}
\begin{split}
&\int_{V_c}|U|^2dV+\int_{V_{c,0,\varepsilon}}|U|^2dV+\int_{ V_{c,1,\varepsilon}}|U|^2dV\\
&\leq(2+8\tilde{C}_1\varepsilon^3+8\tilde{C}_3\varepsilon^3)Q(U,U)
+4C(\varepsilon)\int_{K_\varepsilon}|U|^2dV
+2\int_{K_\varepsilon^\ast}|U|^2dV.
\end{split}
\end{equation}
Let $V_{c,0,1}=V_c\cup V_{c,0,\varepsilon}\cup V_{c,1,\varepsilon}$,
$F=K_\varepsilon\cup K_\varepsilon^\ast$. Applying the approximation
theorem by the smooth forms as mentioned at the end of $\S 3$,  we
thus complete the proof of Theorem \ref{101}. The proof of Corollary
\ref{202} follows from Theorem \ref{101} by the standard argument as
in H\"ormander [Ho1].

\begin{remark}
From the proof of the main theorem,  we actually obtained the
following stronger estimate  containing the boundary term for
$U\in\mathcal{E}_c^{0,q}\cap Dom(T)\cap Dom{S^\ast}$:
\begin{equation}
\int_{\p S^{\e_0}_\varphi}|U|^2dV+\int_{V_{c,0,1}}|U^2dV\leq
C\left(Q(U,U)+\int_F|U|^2dV\right),
\end{equation}
where $F$ is a certain fixed compact subset  of $S^{\e_0}_\varphi$.
 Other related
sub-elliptic estimates  will be discussed in [Li].
\end{remark}
\bigskip
\bibliographystyle{amsalpha}

\noindent Xiaojun Huang (huangx@math.rutgers.edu),
Department of Mathematics, Rutgers University, New Brunswick, NJ
08903, USA.

\medskip
\noindent Xiaoshan Li (xiaoshanli@whu.edu.cn), School of Mathematics
and Statistics, Wuhan University, Hubei 430072,  China.
\end{document}